\documentclass[a4paper,oneside,10pt]{article}%
\usepackage{amsmath}
\usepackage{amsfonts}
\usepackage{amssymb}
\usepackage{graphicx}
\usepackage[square,numbers,sort&compress]{natbib}%
\providecommand{\U}[1]{\protect\rule{.1in}{.1in}}

\newtheorem{theorem}{Theorem}[section]

\newtheorem{definition}[theorem]{Definition}
\newtheorem{assumption}[theorem]{Assumption}

\newtheorem{lemma}[theorem]{Lemma}

\newtheorem{remark}[theorem]{Remark}

\numberwithin{equation}{section}

\begin{document}

\title{Central limit theorems for bounded random variables under belief measures}
\author{Xiaomin Shi\thanks{School of Mathematics and Quantitative Economics,
Shandong University of Finance and Economics, Jinan 250014, China and Zhongtai Securities Institute for Financial Studies, Shandong University, Jinan, 250100, PR China. Email: shixm@mail.sdu.edu.cn. This work was supported by National Natural
Science Foundation of China (No. 11401091).}}
\date{}
\maketitle
\textbf{Abstract}. Recently a new type of central limit theorem for belief functions was given in Epstein et al. \cite{ES1}. In this paper, we generalize the central limit theorem in Epstein et al. \cite{ES1}  to accommodate general bounded random variables.  These results are natural extension of the classical central limit theory for additive probability measures.


\section{Introduction}

Recently, a central limit theorem (CLT for short) for belief functions was given in Epstein et al. \cite{ES1} (Theorem 3.1) to construct suitably robust confidence regions for incomplete models. We state their CLT for the readers' convenience:
\begin{theorem}
\label{EKS}
Let $\Lambda_{\theta_n}\rightarrow \Lambda\in \mathbb{R}^{J \cdot J}$ and $c_n\rightarrow c\in\mathbb{R}^J$. Then
\begin{equation}
\label{cltber}
\nu_{\theta_n}^\infty(\cap_{j=1}^J\{s^\infty:\sqrt{n}[\nu_{\theta_n}(A_j)-\Psi_n(s^\infty)(A_j)]\leq c_{nj}\})\rightarrow \mathbf{N}_J(c;\Lambda).
\end{equation}
\end{theorem}

In the theorem above, structure parameter $\theta$ and $J=1,2......$ be fixed a prior. $S$ is a finite state space, and $A_1, ..., A_J$ are $J$ subsets of $S$. $\Psi_n(s^\infty)(A_j)=\frac{1}{n}\sum\limits_{i=1}^n I_{\{s_i\in A_j\}},\ j=1,...,J$ are the empirical frequency measure of $A_j$ in the first $n$ experiments along the sample $s^\infty=(s_1,s_2,...)$.
\begin{equation}
cov_\theta(A_i,A_j)=\nu_\theta(A_i\cap A_j)-\nu_\theta(A_i)\nu_\theta( A_j).
\end{equation}

$\Lambda_\theta$ is the $J\times J$ symmetric and positive semidefinite matrix $(cov_\theta(A_i,A_j))$ and $\nu_\theta$ is a belief function on $S$.

\[
\mathbf{N}_J(c;\Lambda)=P(\xi\leq c),
\]
where $\xi$ is a $J$-dimensional normal random variable with zero mean and covariance matrix $\Lambda$ and relation $\xi\leq c$ is in the vector sense.

For more clearly, let $ \Lambda_{\theta_n}\equiv \Lambda, c_n\equiv c$, then Theorem \ref{EKS} reduced to
\begin{equation}
\nu^\infty(\cap_{j=1}^J\{s^\infty:\sqrt{n}[\nu(A_j)-\Psi(s^\infty)(A_j)]\leq c\})\rightarrow \mathbf{N}_J(c;\Lambda).
\end{equation}

Suppose further that $J=2,\ A_1\subset S$, $A_2=S\setminus A_1$, then convergence relation (\ref{cltber}) boils down to
\begin{equation}
\label{cltinter}
\nu^\infty(-\frac{c_1}{\sqrt{n}}+\nu(A_1)\leq \Psi_n(s^\infty)(A_1)\leq \nu^*(A_1)+\frac{c_2}{\sqrt{n}})\rightarrow \mathbf{N}_J(c;\Lambda),
\end{equation}
where $c=(c_1,c_2)^\prime$ and
$\Lambda=\left(
\begin{array}
[c]{cc}%
\nu(A_1)(1-\nu(A_1)) & -\nu(A_1)\nu(A_2)\\
-\nu(A_1)\nu(A_2) & \nu(A_2)(1-\nu(A_2))%
\end{array}
\right).  $

Formula (\ref{cltinter}) shows that the belief measure of $\Psi_n(s^\infty)(A_1)$ in a two-sided interval will converge to a bivariate normal distribution which is completely different from the classical CLTs under probabilities.
To our best knowledge, the above theorem is the first result concerning CLTs for two-sided intervals under belief measures. Formula (\ref{cltinter}) dazzles us because we could hardly  imagine the belief measure for a two-sided interval should be approximated by a bivariate normal distribution. Although belief measures, also called infinitely monotone capacities, are the most special non-additive measures (see page 76 in \cite{WKl}), the above result is a milestone in studying CLTs under non-additive probabilities.

It is well known the classical LLNs (laws of large numbers) and CLTs  play a significant role in the
development of probability theory and its applications. The key in the proof
of these theorems lies in the additivity of probabilities. Such additivity assumption is widely used in statistics, economics, decision theory etc, but it is yet not reasonable in many cases. For instance, the famous Allais paradox and  Ellsberg paradox (see Allais \cite{Al} and Ellsberg \cite{El} respectively) in economic decision theory challenged the additivity assumption a lot.
In the seminal paper \cite{Ch}, Choquet proposed the theory of capacities where he removed the additivity assumption. Choquet's theory has been applied to statistics by Huber and Strassen \cite{HS} where Neyman-Pearson lemma for capacities was studied. Also, the Allais paradox and  Ellsberg paradox could be explained partially by capacities, or more generally, non-additive probabilities. As the role of classical limit theories played in the probability theory, limit theories under non-additive probabilities
are urgent for the more potential general applications.

Given a sequence $\{X_{i}\}_{i=1}^{\infty}$ of independent and identically
distributed (i.i.d. for short) random variables for non-additive
probabilities (or capacities), there have been a number of papers related to strong laws of
large numbers. We refer to \cite{CWL},\cite{DW},\cite{ESc}%
,\cite{MM},\cite{Ma},\cite{WF}. But to the best of our knowledge, only a few
papers studied central limit theorems for non-additive probabilities.
Motivated by measuring risk and other financial problems with uncertainty,
Peng \cite{Pe3},\cite{Pe4},\cite{Pe5},\cite{Pe6} put forward the notion of
i.i.d. random variables in the sublinear expectations space ($G$-expectation especially), which accommodated  volatility uncertainty. $G$-expectation had many applications in mathematical finance. For instance, Epstein and Ji \cite{EJ1},\cite{EJ2} formulated a model of utility
that captures the decision-maker's concern with ambiguity about both the
drift and volatility of the driving process under $G$ framework.  And Peng \cite{Pe7},\cite{Pe8} gave a central limit theorem that the accumulated independent and identically distributed random variables converge in law to the $G$-normal distribution. As for CLT and laws of large numbers under sublinear expectation, we also refer to \cite{HC}, \cite{HZ}, \cite{LS}, and \cite{ZC}.

But all the above work consider problems from the perspective of sublinear expectations. As we all known, there is a one-to-one correspondence between the set of linear expectations and the set of  $\sigma$-additive probability measures. But in the sublinear case, this correspondence no longer holds true. Thus CLTs under sublinear expectations neither imply nor are implied by CLTs under non-additive probabilities.
For CLTs under non-additive probabilities, Marinacci \cite{Ma} proves a CLT for a class of capacities that he calls \textquotedblleft controlled". But he only considers CLTs for unidirectional intervals which looks the same as the classical CLTs as long as the partial sums of $\{X_{i}\}_{i=1}^\infty$ is standardized by proper parameters. As we have mentioned,  Theorem 3.1 in Epstein et al. \cite{ES1} is the first result concerning CLTs for two-sided intervals under belief measures. And it has a nice application in Jovanovic entry game and  binary experiments, or more generally, in an environment where the agent does not have enough information to formulate a probability prediction so he can only hold a coarse picture (belief) about the uncertainty. In fact, they construct robust confidence regions for models in which a probabilistic prediction is absent via this central limit theorem.

But Theorem 3.1 in  \cite{ES1} considers only Bernoulli random variables which limit the potential application of the CLTs for belief measures.
If one were in a belief measure world with general (continuum) state space and
general bounded random variables, Theorem \ref{EKS} does not work. So we generalize the CLTs of (\ref{cltinter}) type  to accommodate general bounded random variables. But we don't
involve the concrete applications in this paper. And we think it's still a long way to obtain CLTs under general non-additive probabilities while we hope CLTs under belief measures will give some hints.

Extending from Bernoulli random variables to general bounded random variable is interesting and far from trivial. For instance, if random variables $X$ and $Y$ are of Bernoulli type, the expectations $EX$ and $E(XY)$ are easy to calculate by enumeration. But for general random variables $X$ and $Y$, the expectations $EX$ and $E(XY)$ are not easy to get, especially $E(XY)$.
In the following, random variables $\underline Z_i$ and $\bar Z_i$
are not primitives. Although it is not very hard to translate $E^{P_\nu^\infty}(\underline Z_i)$, $E^{P_\nu^\infty}(\underline Z_i)^2$, $E^{P_\nu^\infty}(\bar Z_i)$ and $E^{P_\nu^\infty}(\bar Z_i)^2$ in terms of the primitive $\nu$. But $E^{P_\nu^\infty}(\underline Z_i\bar Z_i)$ is far from easy to translate.

The paper is organized as follows. In section 2, we give some preliminaries
about belief measures. In section 3, we give a unidirectional central limit
theorem under belief measures and a central limit theorem for two-sided
intervals is given in section 4.

\section{Preliminaries}

Let $\Omega$ be a Polish space and $\mathcal{B}(\Omega)$ be the Borel $\sigma
$-algebra on $\Omega$. We denote by $\mathcal{K}(\Omega)$ the compact subsets
of $\Omega$ and $\Delta(\Omega)$ the space of all probability measures on
$\Omega$. The set $\mathcal{K}(\Omega)$ will be endowed with the Hausdorff
topology generated by the topology of $\Omega$.

\begin{definition}
\label{def2.1} A belief measure on $(\Omega,\mathcal{B}(\Omega))$ is  defined as a set function $\nu: \mathcal{B}(\Omega)\rightarrow[0,1]$ satisfying:

(i) $\nu(\emptyset)=0$ and $\nu(\Omega)=1$;

(ii) $\nu(A)\leq\nu(B)$ for all Borel sets $A\subset B$;

(iii) $\nu(B_{n})\downarrow\nu(B)$ for all sequences of Borel sets
$B_{n}\downarrow B$;

(iv) $\nu(G)=\sup\{\nu(K):K\subset G,K$ is compact$\}$, for all open set G;

(v) $\nu$ is totally monotone (or $\infty$-monotone): for all Borel sets
$B_{1},...,B_{n}$,
\[
\nu(\bigcup_{i=1}^{n}B_{j})\geq\sum_{\emptyset\neq I\subset\{1,...,n\}}%
(-1)^{|I|+1}\nu(\bigcap_{j\in I}B_{j})
\]
where $|I|$ is the number of elements in $I$.
\end{definition}

Its conjugate $\nu^*$ is given by
\[
\nu^*(A)=1-\nu(A^c),
\]
where $A^c$ is the complementary set of $A$.

Note that, by Phillipe et al. \cite{PDJ}, for all $A\in\mathcal{B}(\Omega),\{K\in
\mathcal{K}(\Omega):K\subset A\}$ is universally measurable. Denote by
$\mathcal{B}_{u}(\mathcal{K}(\Omega))$ the $\sigma$-algebra of all subsets of
$\mathcal{K}(\Omega)$ which are universally measurable. The following theorem
belongs to Choquet \cite{Ch}. We also refer to Phillipe et al. \cite{PDJ}.

\begin{theorem}
\label{Choquet} The set function $\nu:\mathcal{B}(\Omega)\rightarrow
\lbrack0,1]$ is a belief measure if and only if there exists a probability
measure $P_{\nu}$ on $(\mathcal{K}(\Omega),\mathcal{B}(\mathcal{K}(\Omega)))$
such that
\[
\nu(A)=P_{\nu}\big(\{K\in\mathcal{K}(\Omega):K\subset A\}\big),~~~\forall
A\in\mathcal{B}(\Omega).
\]
Moreover, there exists a unique extension of $P_{\nu}$ to $(\mathcal{K}%
(\Omega),\mathcal{B}_{u}(\mathcal{K}(\Omega)))$. In the following, we still
denote the extension by $P_{\nu}$.
\end{theorem}

Let $\Omega^\infty$ be the infinite product space of $\Omega$, and $\mathcal{B}(\Omega^\infty)$ the Borel $\sigma$-algebra of $\Omega^\infty$.  For any $A\in\mathcal{B}(\Omega^{\infty})$, we define%
\begin {equation}
\label{nuinfty}
\nu^{\infty}(A)=P_{\nu}^{\infty}\big(\{\widetilde{K}=K_{1}\times K_{2}%
\times...\in(\mathcal{K}(\Omega))^{\infty}:\widetilde{K}\subset A\}\big),
\end{equation}
where $P_{\nu}^{\infty}\in\Delta((\mathcal{K}(\Omega))^{\infty})$ is the
i.i.d. product probability measure. According to Epstein and Seo \cite{ES}
(Lemma A.2.), $\nu^{\infty}$ is the unique belief measure on $(\Omega^{\infty
},\mathcal{B}(\Omega^{\infty}))$ corresponding to $P_{\nu}^{\infty}$. The conjugate of $\nu^\infty$ is denoted by $\nu^{*\infty}$.

\begin{definition}
\label{defidentical} The elements of a sequence $\{\xi_{i}\}_{i=1}^{\infty}$
of random variables on $(\Omega,\mathcal{B}(\Omega))$ are called identically
distributed w.r.t. $\nu$, denoted by $\xi_{i}\overset{d}{=}\xi_{j}$, if for
each $i,j\geq1$ and for all intervals $G$ of $\mathbb{R}$,
\[
\nu(\xi_{i}\in G)=\nu(\xi_{j}\in G).
\]

\end{definition}
\begin{remark}
Note that  belief functions are NOT  additive in general, so they are NOT characterized by their values on intervals only. But we still use \textquotedblleft all intervals\textquotedblright\  rather than \textquotedblleft all Borel sets\textquotedblright\  in Definition 2.3 because our CLTs deal with intervals only.
\end{remark}

In this paper, we are mainly interested in central limit theorems for  random variables in the following form.

\begin{assumption}
\label{assu-1}
Let $\{Y_{i}\}_{i=1}^\infty$ be a sequence of bounded random variables on $(\Omega,\mathcal{B}(\Omega))$, $Y_{i}\overset{d}{=}%
Y_{j},i, j\geq1$. We denote $X_{i}(\omega_{1},\omega_{2},...):=Y_{i}(\omega
_{i}),\ i\geq1,\ (\omega_{1},\omega_{2},...)\in \Omega^\infty$. And let $M$ be the least upper bound of $\{Y_{i}\}_{i=1}^\infty$, that is, $M=inf\{C \geq 0:\nu(|Y_i|\leq C)=1\}$.
\end{assumption}

Under Assumption \ref{assu-1}, $\{X_{i}\}_{i=1}^\infty$ is a sequence of random variables on $(\Omega^{\infty},\mathcal{B}(\Omega^{\infty}))$ bounded by $M$.

It's well known that the Choquet integral of $Y_i, i\geq 1$ w.r.t. belief measure $\nu$, denoted by $\mathcal{E}_\nu(Y_i)$ is
\begin{align}
\label{chointe}
\mathcal{E}_\nu(Y_i)&=\int_{0}^{M} \nu(Y_{i}\geq t)\mathrm{d}t+\int_{-M}^{0} [\nu(Y_{i}\geq t)-1]\mathrm{d}t \nonumber \\
&=\int_{0}^{M} \nu(Y_{i}> t)\mathrm{d}t+\int_{-M}^{0} [\nu(Y_{i}> t)-1]\mathrm{d}t, i\geq 1.
\end{align}

In the remaining paper, we shall use equation (\ref{chointe}) repeatedly without reference.

\section{Unidirectional central limit theorem}

This section mainly studies a central limit theorem for belief measure in a unidirectional interval. As we will show, the belief measures of the events where the partial sums of $\{X_{i}\}_{i=1}^\infty$ is in  unidirectional intervals can be approximated by a normal distribution in another probability space. But
the belief measures
of the events where the partial sums of $\{X_{i}\}_{i=1}^\infty$ is in   two-sided intervals
 are rather involved and can not approximated by a normal distribution. Actually it should be approximated by a bivariate normal distribution. So we leave the central limit theorem for a two-sided interval to the next section.

\begin{lemma}
\label{nuinf}
Let $\{Y_{i}\}_{i=1}^\infty$ and $\{X_{i}\}_{i=1}^\infty$  be defined in Assumption \ref{assu-1}, then for any interval $G$ of $\mathbb{R}$, we have $$\nu^\infty(X_i\in G)=\nu(Y_i\in G),\ i\geq 1.$$
\end{lemma}
\noindent\textbf{Proof. }
In fact, by (\ref{nuinfty}) and Theorem \ref{Choquet},
\begin{align*}
\nu^\infty(X_i\in G)&=P_{\nu}^{\infty}\big(\{\widetilde{K}=K_{1}\times K_{2}%
\times...\in(\mathcal{K}(\Omega))^{\infty}:\widetilde{K}\subset \{(\omega_1, \omega_2, ...): X_i(\omega_1, \omega_2, ...)\in G\}\}\big)\\
&=P_{\nu}^{\infty}\big(\{\widetilde{K}=K_{1}\times K_{2}%
\times...\in(\mathcal{K}(\Omega))^{\infty}:\widetilde{K}\subset \{(\omega_1, \omega_2, ...): Y_i(\omega_i)\in G\}\}\big)\\
&=P_{\nu}^{\infty}\big(\{\widetilde{K}=K_{1}\times K_{2}%
\times...\in(\mathcal{K}(\Omega))^{\infty}: K_i\subset \{\omega_i: Y_i(\omega_i)\in G\}\}\big)\\
&=P_{\nu}\big(\{K_{i}\in\mathcal{K}(\Omega): K_i\subset \{\omega_i: Y_i(\omega_i)\in G\}\}\big)\\
&=\nu(Y_i\in G),\ \  i\geq 1.  \ \ \ \ \ \Box
\end{align*}

\begin{lemma}
\label{lemma1}
For $\{Y_{i}\}_{i=1}^\infty$ and $\{X_{i}\}_{i=1}^\infty$ defined in Assumption \ref{assu-1}, let
\begin{align*}
&  \underline Z_{i}(K_{1}\times K_{2}\times...)=\inf\limits_{\omega_{i}\in K_{i}}%
X_{i}(\omega_{1},\omega_{2},...),\\
&  \bar{Z}_{i}(K_{1}\times K_{2}\times...)=\sup\limits_{\omega_{i}\in K_{i}%
}X_{i}(\omega_{1},\omega_{2},...), \ K_{i}\in\mathcal{K}(\Omega),~~i\geq1.
\end{align*}
We claim that $\{\underline Z_{i}\}_{i=1}^{\infty}$ and $\{\bar{Z}_{i}\}_{i=1}^{\infty}$
are two sequences of i.i.d. bounded random variables defined on $\big((\mathcal{K}%
(\Omega))^{\infty},\mathcal{B}_{u}((\mathcal{K}(\Omega))^{\infty})\big)$
w.r.t. $P_{\nu}^{\infty}$. Moreover,
\begin{equation}
E^{P_{\nu}^{\infty}}[\underline Z_{i}]=\mathcal{E}_{\nu}(Y_i)=:\underline{\mu}, \ E^{P_{\nu}^{\infty}}[\bar{Z}_{i}]=\mathcal{E}_{\nu^{*}}(Y_i)=:\bar{\mu}\;,~~i\geq1, \label{ExpVar}%
\end{equation}%
\begin{equation}
Var^{P_{\nu}^{\infty}}%
[\underline Z_{i}]=\int_{0}^{M}2t\nu(Y_{i}\geq t)\mathrm{d}t+\int_{-M}^{0}%
2t[\nu(Y_{i}\geq t)-1]\mathrm{d}t-\underline{\mu}^{2}=:\underline\sigma^2,
\end{equation}
\begin{equation}
Var^{P_{\nu}^{\infty}}[\bar
{Z}_{i}]=\int_{0}^{M}2t\nu^{*}(Y_{i}\geq t)\mathrm{d}t+\int_{-M}^{0}%
2t[\nu^{*}(Y_{i}\geq t)-1]\mathrm{d}t-\underline{\mu}^{2}=:\bar{\sigma}^{2},~~i\geq1. \label{ExpVar-2}%
\end{equation}
\end{lemma}

\noindent\textbf{Proof. }
We only prove the assertions for $\{\underline Z_{i}\}_{i=1}^{\infty}$. The proofs for
$\{\bar{Z}_{i}\}_{i=1}^{\infty}$ are parallel.

By Assumption \ref{assu-1}, $M$ is the least upper bound of $\{X_i\}_{i=1}^\infty$ in the sense that
\[
M=inf\{C \geq 0:\nu^\infty(|X_i|\leq C)=1\}, \ i\geq 1.
\]
So $\{\underline Z_{i}\}_{i=1}^{\infty}$ is bounded by $M$.

For each $i,j\geq1,\ \forall t,t_{i},t_{j}\in\mathbb{R}$,
\begin{align*}
P_{\nu}^{\infty}(\underline Z_{i}\geq t)  &  =P_{\nu}^{\infty}\big(\{\widetilde{K}%
=K_{1}\times K_{2}\times...\in(\mathcal{K}(\Omega))^{\infty}:\underline Z_{i}%
(\widetilde{K})\geq t\}\big)\\
&  =P_{\nu}^{\infty}\big(\{\widetilde{K}=K_{1}\times K_{2}\times
...\in(\mathcal{K}(\Omega))^{\infty}:\inf_{\omega_{i}\in K_{i}}X_{i}%
(\omega_{1},\omega_{2},...)\geq t\}\big)\\
&  =P_{\nu}^{\infty}\big(\{\widetilde{K}=K_{1}\times K_{2}\times
...\in(\mathcal{K}(\Omega))^{\infty}:\widetilde{K}\subset\{(\omega_{1}%
,\omega_{2},...):X_{i}(\omega_{1},\omega_{2},...)\geq t\}\}\big)\\
&  =\nu^{\infty}(X_{i}\geq t)\\
&  =\nu(Y_{i}\geq t)\\
&  =\nu(Y_{j}\geq t)\\
&  =\nu^{\infty}(X_{j}\geq t)\\
&  =P_{\nu}^{\infty}(\underline Z_{j}\geq t),
\end{align*}
and for $i\neq j,$%
\[%
\begin{array}
[c]{rl}
& P_{\nu}^{\infty}(\underline Z_{i}\geq t_{i},\underline Z_{j}\geq t_{j})\\
= & P_{\nu}^{\infty}\big(\{\widetilde{K}=K_{1}\times K_{2}\times
...\in(\mathcal{K}(\Omega))^{\infty}:\underline Z_{i}(\widetilde{K})\geq t_{i}%
,\underline Z_{j}(\widetilde{K})\geq t_{j}\}\big)\\
= & P_{\nu}^{\infty}\big(\{\widetilde{K}=K_{1}\times K_{2}\times
...\in(\mathcal{K}(\Omega))^{\infty}:\inf_{\omega_{i}\in K_{i}}X_{i}%
(\omega_{1},\omega_{2},...)\geq t_{i},\inf_{\omega_{j}\in K_{j}}X_{j}%
(\omega_{1},\omega_{2},...)\geq t_{j}\}\big)\\
= & P_{\nu}^{\infty}\big(\{\widetilde{K}=K_{1}\times K_{2}\times
...\in(\mathcal{K}(\Omega))^{\infty}:\inf_{\omega_{i}\in K_{i}}X_{i}%
(\omega_{1},\omega_{2},...)\geq t_{i}\}\big)\times\\
& P_{\nu}^{\infty}\big(\{\widetilde{K}=K_{1}\times K_{2}\times...\in
(\mathcal{K}(\Omega))^{\infty}:\inf_{\omega_{j}\in K_{j}}X_{j}(\omega
_{1},\omega_{2},...)\geq t_{j}\}\big)\\
= & P_{\nu}^{\infty}(\underline Z_{i}\geq t_{i})P_{\nu}^{\infty}(\underline Z_{j}\geq t_{j}),
\end{array}
\]
where the third equality follows from the definition of product  measure.

Similarly, we can easily get for any different positive integer $i_i,...,i_n$,
\[
P_{\nu}^{\infty}(\underline Z_{i_1}\geq t_{i_1},...,\underline Z_{i_n}\geq t_{i_n})=P_{\nu}^{\infty}(\underline Z_{i_1}\geq t_{i_1})\times \cdot\cdot\cdot \times P_{\nu}^{\infty}(\underline Z_{i_n}\geq t_{i_n}).
\]
Thus, $\{\underline Z_{i}\}_{i=1}^{\infty}$ is a sequence of  independent and
identically distributed random variables under $P_{\nu}^{\infty}$.

Moreover, we calculate the expectation and variation of $\underline Z_{i}$ as follows:
\begin{align*}
E^{P_{\nu}^{\infty}}(\underline Z_{i})  &  =\int_{0}^{M} P_{\nu}^{\infty}(\underline Z_{i}\geq
t)\mathrm{d}t+\int_{-M}^{0} [P_{\nu}^{\infty}(\underline Z_{i}\geq t)-1]\mathrm{d}t \ \\
&  =\int_{0}^{M} \nu^{\infty}(X_{i}\geq t)\mathrm{d}t+\int_{-M}^{0} [\nu^{\infty}(X_{i}\geq t)-1]\mathrm{d}t \ \\
&  =\int_{0}^{M} \nu(Y_{i}\geq t)\mathrm{d}t+\int_{-M}^{0} [\nu(Y_{i}\geq t)-1]\mathrm{d}t \ \\
&  =\mathcal{E}_{\nu}(Y_{i}), ~~i\geq1,
\end{align*}
where the third equality is due to Lemma \ref{nuinf}.
And
\begin{align*}
Var^{P_{\nu}^{\infty}}(\underline Z_{i})  &  =E^{P_{\nu}^{\infty}}(\underline Z_{i}^{2}%
)-\big(E^{P_{\nu}^{\infty}}(\underline Z_{i})\big)^{2}\\
&  =\int_{0}^{M^{2}}P_{\nu}^{\infty}(\underline Z_{i}^{2}\geq t)\mathrm{d}%
t-\underline{\mu}^{2}\\
&  =\int_{0}^{M^{2}}P_{\nu}^{\infty}(\underline Z_{i}\geq\sqrt{t})\mathrm{d}t+\int%
_{0}^{M^{2}}P_{\nu}^{\infty}(\underline Z_{i}\leq-\sqrt{t})\mathrm{d}t-\underline{\mu
}^{2}\\
&  =\int_{0}^{M^{2}}\nu^{\infty}(X_{i}\geq\sqrt{t})\mathrm{d}t+\int_{0}%
^{M^{2}}[1-\nu^{\infty}(X_{i}\geq-\sqrt{t})]\mathrm{d}t-\underline{\mu}^{2}\\
&  =\int_{0}^{M}2t\nu^{\infty}(X_{i}\geq t)\mathrm{d}t+\int_{-M}^{0}%
2t[\nu^{\infty}(X_{i}\geq t)-1]\mathrm{d}t-\underline{\mu}^{2} \\
&  =\int_{0}^{M}2t\nu(Y_{i}\geq t)\mathrm{d}t+\int_{-M}^{0}%
2t[\nu(Y_{i}\geq t)-1]\mathrm{d}t-\underline{\mu}^{2},~~i\geq1,
\end{align*}
where the last equality is due to Lemma \ref{nuinf}. $\ \ \ \ \ \Box$
\begin{remark}
In probability theory, $E(Y_i^2)-(EY_i)^2$ is a variant  definition of variance. Under non-additive probabilities, there are many ways to define the \textquotedblleft variance" among which $\mathcal{E}_{\nu}(Y_i^2)-(\mathcal{E}_{\nu} (Y_i))^2$ seems intuitive.
It's easy to verify that for $0\leq Y_i\leq M,$
\begin{align*}
\underline\sigma^2&=\int_{0}^{M}2t\nu(Y_{i}\geq
t)\mathrm{d}t+\int_{-M}^{0}2t[\nu(Y_{i}\geq t)-1]\mathrm{d}%
t-\underline{\mu}^{2}\\
&=\mathcal{E}_{\nu}(Y_i^2)-(\mathcal{E}_{\nu} (Y_i))^2.
\end{align*}
Unfortunately $\underline\sigma^2\neq\mathcal{E}_{\nu}(Y_i^2)-(\mathcal{E}_{\nu} (Y_i))^2$ for general $|Y_i|\leq M$.
For example, readers can easily verify that for $-M\leq Y_i\leq 0$,
\begin{align*}
\underline\sigma^2&=\int_{0}^{M}2t\nu(Y_{i}\geq
t)\mathrm{d}t+\int_{-M}^{0}2t[\nu(Y_{i}\geq t)-1]\mathrm{d}%
t-\underline{\mu}^{2}\\
&=\mathcal{E}_{\nu^{*}}(Y_i^2)-(\mathcal{E}_{\nu} (Y_i))^2\neq \mathcal{E}_{\nu}(Y_i^2)-(\mathcal{E}_{\nu} (Y_i))^2,
\end{align*}
except that $\nu$ degenerates to a probability measure.
\end{remark}

Let $\mathbf{N}(\cdot)$ be the cumulate distribution function of a standard normal distribution. Our main result in this section is the following central limit theorem.

\begin{theorem}
\label{Thm-clt} For $\{X_{i}\}_{i=1}^\infty$ defined in Assumption \ref{assu-1}, we have
\begin{equation}
\lim_{n\rightarrow\infty}\nu^{\infty}\Bigg(\frac{\sum\limits_{i=1}^{n}%
X_{i}-n\underline{\mu}}{\sqrt{n}\underline{\sigma}}\geq\alpha
\Bigg)=1-\mathbf{N}(\alpha), \label{clt1}%
\end{equation}
and
\begin{equation}
\lim_{n\rightarrow\infty}\nu^{\infty}\Bigg(\frac{\sum\limits_{i=1}^{n}%
X_{i}-n\bar{\mu}}{\sqrt{n}\bar{\sigma}}<\alpha\Bigg)=\mathbf{N}(\alpha
),~~\forall\alpha\in\mathbb{R}, \label{clt2}%
\end{equation}
where $\underline\mu, \ \bar\mu, \ \underline\sigma, \ \bar\sigma$ are defined in Lemma \ref{lemma1}.
\end{theorem}

\noindent\textbf{Proof. }
By Theorem \ref{Choquet}, we have
\begin{align*}
& \ \ \ \ \  \nu^{\infty}\Big(\sum\limits_{i=1}^{n} X_{i}\geq\alpha\Big)\\
&  =P_{\nu}^{\infty}\big(\{\widetilde{K}=K_{1}\times K_{2}\times
...\in(\mathcal{K}(\Omega))^{\infty}:\widetilde{K}\in\{(\omega_{1},\omega
_{2},...):\sum_{i=1}^{n} X_{i}\geq\alpha\}\}\big)\\
&  =P_{\nu}^{\infty}\big(\{\widetilde{K}=K_{1}\times K_{2}\times
...\in(\mathcal{K}(\Omega))^{\infty}:\widetilde{K}\in\{(\omega_{1},\omega
_{2},...):\sum_{i=1}^{n} Y_{i}(\omega_{i})\geq\alpha\}\}\big)\\
&  =P_{\nu}^{\infty}\big(\{\widetilde{K}=K_{1}\times K_{2}\times
...\in(\mathcal{K}(\Omega))^{\infty}:\inf_{(\omega_{1},\omega_{2},...)\in
K_{1}\times K_{2}\times...}\sum_{i=1}^{n} Y_{i}(\omega_{i})\geq\alpha\}\big)\\
&  =P_{\nu}^{\infty}\big(\{\widetilde{K}=K_{1}\times K_{2}\times
...\in(\mathcal{K}(\Omega))^{\infty}:\sum_{i=1}^{n} \inf_{\omega_{i}\in K_{i}}
Y_{i}(\omega_{i})\geq\alpha\}\big)\\
&  =P_{\nu}^{\infty}\big(\{\widetilde{K}=K_{1}\times K_{2}\times
...\in(\mathcal{K}(\Omega))^{\infty}:\sum_{i=1}^{n} \underline Z_{i}(K_{1}\times
K_{2}\times...)\geq\alpha\}\big)\\
&  =P_{\nu}^{\infty}\Big(\sum_{i=1}^{n} \underline Z_{i}\geq\alpha\Big), ~~\forall
\alpha\in\mathbb{R}.
\end{align*}

Thus
\begin{align}
\nu^{\infty}\Bigg(\frac{\sum\limits_{i=1}^{n}X_{i}-n\underline{\mu}}{\sqrt
{n}\underline{\sigma}}\geq\alpha\Bigg)=P_{\nu}^{\infty}\Bigg(\frac
{\sum\limits_{i=1}^{n}\underline Z_{i}-n\underline{\mu}}{\sqrt{n}\underline{\sigma}}%
\geq\alpha\Bigg), ~~\forall\alpha\in\mathbb{R}.
\end{align}

Applying the classical central limit theorem to $\{\underline Z_{i}\}_{i=1}^{\infty}$
deliver the corresponding limit theorem (\ref{clt1}).

As for the proof of limit theorem (\ref{clt2}), we just reverse the inequality
sign and replace $\{\underline Z_{i}\}_{i=1}^{\infty}$ by $\{\bar{Z}_{i}\}_{i=1}^{\infty
}$ in the preceding proof. $\ \ \ \ \ \Box$

\section{Central limit theorem for two-sided intervals}

In this section we will give the  CLT for belief measures in two-sided intervals. We approximate the belief measures
of the events where the partial sums of $\{X_{i}\}_{i=1}^\infty$ is in a two-sided intervals
by a bivariate normal distribution.

Let $\mathbf{N}_{2}(\cdot,\cdot;\rho)$ be the cumulate distribution function
for the bivariate normal distribution with zero means, unit variances and correlation
coefficient $\rho$, i.e.,
\[
\mathbf{N}_{2}(\alpha_{1},\alpha_{2};\rho)=\mathrm{Pr}(B_{1}%
\leq\alpha_{1},B_{2}\leq\alpha_{2}),
\]
where $(B_{1},B_{2})$ is a bivariate normal random variable with the
indicated moments.

\begin{lemma}
\label{lemma2}
Under the assumption of Lemma \ref{lemma1}, we have
\begin{align}
\label{rho}
\rho:=corr^{P_{\nu}^{\infty}}(\underline Z_{i},\bar{Z}_{i})=\frac
{M^{2}-M\bar{\mu}+M\underline{\mu}-\int_{-M}^{M}\int_{-M}^{t_{2}}\nu(t_{1}\leq Y_{i}\leq t_{2})\mathrm{d}%
t_{1}\mathrm{d}t_{2}-\underline{\mu}\bar{\mu}%
}{\underline{\sigma}\bar{\sigma}}.
\end{align}
\end{lemma}
\noindent\textbf{Proof. }Note that for positive random variables $\eta$ and $\zeta$,
\begin{align*}
E[\eta\zeta] &=E\Big[\int_0^\infty I_{\{\eta\geq x\}}\mathrm{d}x\int_0^\infty I_{\{\zeta\geq y\}}\mathrm{d}y\Big]\\
&=\int_0^\infty\int_0^\infty P(\eta\geq x, \zeta\geq y)\mathrm{d}x\mathrm{d}y,
\end{align*}
where the second equality follows from Fubini's theorem.
Then
\[
E^{P_{\nu}^{\infty}}[\underline Z_i \bar Z_i]=E^{P_{\nu}^{\infty}}[\underline Z_i^{+} \bar Z_i^{+}]-E^{P_{\nu}^{\infty}}[\underline Z_i^{+} \bar Z_i^{-}]-E^{P_{\nu}^{\infty}}[\underline Z_i^{-} \bar Z_i^{+}]+E^{P_{\nu}^{\infty}}[\underline Z_i^{-} \bar Z_i^{-}]=I_1+I_2+I_3+I_4,
\]
where
\begin{align*}
I_1&=E^{P_{\nu}^{\infty}}[\underline Z_i^{+} \bar Z_i^{+}]\\
&=\int_0^M\int_0^M P_\nu^\infty(\underline Z_i\geq t_1, \bar Z_i\geq t_2)\mathrm{d}t_1\mathrm{d}t_2\\
&=\int_0^M\int_0^M P_\nu^\infty(\underline Z_i\geq t_1)\mathrm{d}t_1\mathrm{d}t_2-\int_0^M\int_0^M P_\nu^\infty(\underline Z_i\geq t_1, \bar Z_i\leq t_2)\mathrm{d}t_1\mathrm{d}t_2,\\
I_2&=-E^{P_{\nu}^{\infty}}[\underline Z_i^{+} \bar Z_i^{-}]\\
&=-\int_0^M\int_0^M P_\nu^\infty(\underline Z_i\geq t_1, -\bar Z_i\geq t_2)\mathrm{d}t_1\mathrm{d}t_2\\
&=-\int_{-M}^0\int_0^M P_\nu^\infty(\underline Z_i\geq t_1, \bar Z_i\leq t_2)\mathrm{d}t_1\mathrm{d}t_2,\\
I_3&=-E^{P_{\nu}^{\infty}}[\underline Z_i^{-} \bar Z_i^{+}]\\
&=-\int_0^M\int_0^M P_\nu^\infty(-\underline Z_i\geq t_1, \bar Z_i\geq t_2)\mathrm{d}t_1\mathrm{d}t_2\\
&=-\int_0^M\int_{-M}^0 P_\nu^\infty(\underline Z_i\leq t_1, \bar Z_i\geq t_2)\mathrm{d}t_1\mathrm{d}t_2\\
&=\int_0^M\int_{-M}^0 \Big[P_\nu^\infty(\bar Z_i\leq t_2)- P_\nu^\infty(\underline Z_i\leq t_1)- P_\nu^\infty(\underline Z_i\geq t_1, \bar Z_i\leq t_2)\Big]\mathrm{d}t_1\mathrm{d}t_2,\\
I_4&=E^{P_{\nu}^{\infty}}[\underline Z_i^{-} \bar Z_i^{-}]\\
&=\int_0^M\int_0^M P_\nu^\infty(-\underline Z_i\geq t_1, -\bar Z_i\geq t_2)\mathrm{d}t_1\mathrm{d}t_2\\
&=\int_{-M}^0\int_{-M}^0 P_\nu^\infty(\underline Z_i\leq t_1, \bar Z_i\leq t_2)\mathrm{d}t_1\mathrm{d}t_2\\
&=\int_{-M}^0\int_{-M}^0 \Big[P_\nu^\infty( \bar Z_i\leq t_2)-P_\nu^\infty(\underline Z_i\geq t_1, \bar Z_i\leq t_2)\Big]\mathrm{d}t_1\mathrm{d}t_2.
\end{align*}

Adding $I_1,I_2,I_3,I_4$ up, we get
\begin{align*}
& \ \ \ \ \ E^{P_{\nu}^{\infty}}[\underline Z_i \bar Z_i]\\
&=-\int_{-M}^M\int_{-M}^M P_\nu^\infty(\underline Z_i\geq t_1, \bar Z_i\leq t_2)\mathrm{d}t_1\mathrm{d}t_2+\int_0^M\int_0^M P_\nu^\infty(\underline Z_i\geq t_1)\mathrm{d}t_1\mathrm{d}t_2\\
& \ \ \ \ \ +\int_0^M\int_{-M}^0 \Big[P_\nu^\infty(\bar Z_i\leq t_2)- P_\nu^\infty(\underline Z_i\leq t_1)\mathrm{d}t_1\mathrm{d}t_2
+\int_{-M}^0\int_{-M}^0 P_\nu^\infty( \bar Z_i\leq t_2)\mathrm{d}t_1\mathrm{d}t_2\\
&  =-\int_{-M}^{M}\int_{-M}^{M}P_{\nu}^{\infty}(t_{1}\leq \underline Z_{i}\leq\bar{Z}%
_{i}\leq t_{2})\mathrm{d}t_{1}\mathrm{d}t_{2}+M\int_{-M}^{M}P_{\nu}^{\infty
}(\bar{Z}_{i}\leq t_{2})\mathrm{d}t_{2}-M\int_{-M}^{M}P_{\nu}^{\infty}%
(\underline Z_{i}\leq t_{1})\mathrm{d}t_{1}+M^{2}\\
&  =-\int_{-M}^{M}\int_{-M}^{M}\nu^{\infty}(t_{1}\leq X_{i}\leq t_{2}%
)\mathrm{d}t_{1}\mathrm{d}t_{2}+M\int_{-M}^{M}P_{\nu}^{\infty}(\bar{Z}_{i}\leq
t_{2})\mathrm{d}t_{2}-M\int_{-M}^{M}P_{\nu}^{\infty}(\underline Z_{i}\leq t_{1}%
)\mathrm{d}t_{1}+M^{2}\\
&  =-\int_{-M}^{M}\int_{-M}^{M}\nu(t_{1}\leq Y_{i}\leq t_{2})\mathrm{d}%
t_{1}\mathrm{d}t_{2}+M(M-\bar{\mu})-M(M-\underline{\mu})+M^{2}\\
&  =M^{2}-M\bar{\mu}+M\underline{\mu}-\int_{-M}^{M}\int_{-M}^{t_{2}}\nu
(t_{1}\leq Y_{i}\leq t_{2})\mathrm{d}t_{1}\mathrm{d}t_{2}.
\end{align*}

It yields that
\begin{align*}
\rho:=corr^{P_{\nu}^{\infty}}(\underline Z_{i},\bar{Z}_{i})&=\frac{E^{P_{\nu}^{\infty}}(\underline Z_i\bar Z_i)-E^{P_{\nu}^{\infty}}(\underline Z_i)E^{P_{\nu}^{\infty}}(\bar Z_i)}{\sqrt{E^{P_{\nu}^{\infty}}(\underline Z_i^2)-(E^{P_{\nu}^{\infty}}(\underline Z_i))^2}\sqrt{E^{P_{\nu}^{\infty}}(\overline Z_i^2)-(E^{P_{\nu}^{\infty}}(\overline Z_i))^2}}\\
&=\frac
{M^{2}-M\bar{\mu}+M\underline{\mu}-\int_{-M}^{M}\int_{-M}^{t_{2}}\nu(t_{1}\leq Y_{i}\leq t_{2})\mathrm{d}%
t_{1}\mathrm{d}t_{2}-\underline{\mu}\bar{\mu}%
}{\underline{\sigma}\bar{\sigma}}. \ \ \ \ \ \Box
\end{align*}

\begin{theorem}
\label{twosideclt} For $\{X_{i}\}_{i=1}^\infty$ defined in Assumption \ref{assu-1}, there is a constant K which does not
depend on $\alpha_{1},\alpha_{2}$ or n, such that,
\[
\Big|\nu^{\infty}(\alpha_{1}\sqrt{n}\underline{\sigma}+n\underline{\mu}%
\leq\sum\limits_{i=1}^{n}X_{i}\leq\alpha_{2}\sqrt{n}\bar{\sigma}+n\bar{\mu
})-\mathbf{N_{2}}(-\alpha_{1},\alpha_{2};-\rho)\Big|\leq\frac{K}{\sqrt{n}},
\]
where $\rho$ is defined in Lemma \ref{lemma2}.

Moreover, the similar result holds if $\alpha_{1}$ and
$\alpha_{2}$ depend on n.
\end{theorem}

\noindent\textbf{Proof: } Let $\{\underline Z_{i}\}_{i=1}^{\infty}$ and $\{\bar{Z}%
_{i}\}_{i=1}^{\infty}$  be the same as in Lemma \ref{lemma1}, and $\rho$ be defined in Lemma \ref{lemma2}.
Then we have
\[
E^{P_{\nu}^{\infty}}\left(
\begin{array}
[c]{cc}%
\frac{\underline Z_{i}-\underline{\mu}}{\underline{\sigma}} & \\
\frac{\bar{Z}_{i}-\bar{\mu}}{\bar{\sigma}} &
\end{array}
\right)  =0,\text{ \ \ }Var^{P_{\nu}^{\infty}}\left(
\begin{array}
[c]{cc}%
\frac{\underline Z_{i}-\underline{\mu}}{\underline{\sigma}} & \\
\frac{\bar{Z}_{i}-\bar{\mu}}{\bar{\sigma}} &
\end{array}
\right)  =\left(
\begin{array}
[c]{cc}%
1 & \rho\\
\rho & 1
\end{array}
\right).  \ \ \
\]

And we have,
\begin{align}
\label{Z} & \ \ \ \ \  \nu^{\infty}(\alpha_{1}\sqrt{n}\underline{\sigma}+n\underline{\mu
}\leq\sum\limits_{i=1}^{n} X_{i}\leq\alpha_{2}\sqrt{n}\bar\sigma+n\bar\mu)
\ \nonumber\\
&  =P_{\nu}^{\infty}\big(\{\widetilde{K}=K_{1}\times K_{2}\times
...\in(\mathcal{K}(\Omega))^{\infty}:\widetilde{K}\subset\{(\omega_{1}%
,\omega_{2},...):\alpha_{1}\sqrt{n}\underline{\sigma}+n\underline{\mu}\leq
\sum\limits_{i=1}^{n} X_{i}\leq\alpha_{2}\sqrt{n}\bar\sigma+n\bar
\mu\}\}\big) \ \nonumber\\
&  =P_{\nu}^{\infty}\big(\{\widetilde{K}=K_{1}\times K_{2}\times
...\in(\mathcal{K}(\Omega))^{\infty}:\widetilde{K}\subset\{(\omega_{1}%
,\omega_{2},...):\alpha_{1}\sqrt{n}\underline{\sigma}+n\underline{\mu}\leq
\sum\limits_{i=1}^{n} Y_{i}(\omega_{i})\leq\alpha_{2}\sqrt{n}\bar\sigma
+n\bar\mu\}\}\big) \ \nonumber\\
&  =P_{\nu}^{\infty}\big(\{\widetilde{K}=K_{1}\times K_{2}\times
...\in(\mathcal{K}(\Omega))^{\infty}:\alpha_{1}\sqrt{n}\underline{\sigma
}+n\underline{\mu}\leq\inf_{(\omega_{1},\omega_{2},...)\in K_{1}\times
K_{2}\times...}\sum\limits_{i=1}^{n} Y_{i}(\omega_{i}) \ \nonumber\\
&  ~~ ~~ ~~\leq\sup_{(\omega_{1},\omega_{2},...)\in K_{1}\times K_{2}%
\times...}\sum\limits_{i=1}^{n} Y_{i}(\omega_{i})\leq\alpha_{2}\sqrt{n}%
\bar\sigma+n\bar\mu\}\big) \ \nonumber\\
&  =P_{\nu}^{\infty}\big(\{\widetilde{K}=K_{1}\times K_{2}\times
...\in(\mathcal{K}(\Omega))^{\infty}:\alpha_{1}\sqrt{n}\underline{\sigma
}+n\underline{\mu}\leq\sum\limits_{i=1}^{n} \inf_{\omega_{i}\in K_{i}}
Y_{i}(\omega_{i})\leq\sum\limits_{i=1}^{n} \sup_{\omega_{i}\in K_{i}}
Y_{i}(\omega_{i})\leq\alpha_{2}\sqrt{n}\bar\sigma+n\bar\mu
\}\big) \ \nonumber\\
&  =P_{\nu}^{\infty}\big(\{\widetilde{K}=K_{1}\times K_{2}\times
...\in(\mathcal{K}(\Omega))^{\infty}:\alpha_{1}\sqrt{n}\underline{\sigma
}+n\underline{\mu}\leq\sum\limits_{i=1}^{n} \underline Z_{i}(\widetilde{K})\leq
\sum\limits_{i=1}^{n} \bar Z_{i}(\widetilde{K})\leq\alpha_{2}\sqrt{n}%
\bar\sigma+n\bar\mu\}\big) \ \nonumber\\
&  =P_{\nu}^{\infty}(\alpha_{1}\sqrt{n}\underline{\sigma}+n\underline{\mu}%
\leq\sum\limits_{i=1}^{n} \underline Z_{i}(\widetilde{K})\leq\sum\limits_{i=1}^{n} \bar
Z_{i}(\widetilde{K})\leq\alpha_{2}\sqrt{n}\bar\sigma+n\bar\mu) \ \nonumber\\
&  =P_{\nu}^{\infty}\Bigg(\alpha_{1}\leq\frac{\sum\limits_{i=1}^{n}
\underline Z_{i}-n\underline{\mu}}{\sqrt{n}\underline{\sigma}},\frac{\sum\limits_{i=1}%
^{n} \bar Z_{i}-n\bar\mu}{\sqrt{n}\bar\sigma}\leq\alpha_{2}\Bigg).
\end{align}

Using similar method as  in the proof of Lemma \ref{lemma1}, we can easily get $\{(\underline Z_i, \bar Z_i)^{'}\}_{i=1}^\infty$ is a sequence of  i.i.d. random vectors under $P_\nu^\infty$.
Thus, the classical central limit theorem goes through.

Let $T=\left(
\begin{array}
[c]{cc}%
1 & 0\\
-\frac{\rho}{\sqrt{1-\rho^{2}}} & \frac{1}{\sqrt{1-\rho^{2}}}%
\end{array}
\right)  $. Then
\[
E^{P_{\nu}^{\infty}}\Bigg[T\left(
\begin{array}
[c]{c}%
\frac{\underline Z_{i}-\underline{\mu}}{\underline{\sigma}}\\
\frac{\bar{Z}_{i}-\bar{\mu}}{\bar{\sigma}}%
\end{array}
\right)  \Bigg]=0,\text{ \ \ }Var^{P_{\nu}^{\infty}}\Bigg[T\left(
\begin{array}
[c]{c}%
\frac{\underline Z_{i}-\underline{\mu}}{\underline{\sigma}}\\
\frac{\bar{Z}_{i}-\bar{\mu}}{\bar{\sigma}}%
\end{array}
\right)  \Bigg]=\left(
\begin{array}
[c]{cc}%
1 & 0\\
0 & 1
\end{array}
\right)  ,\ \ \
\]
and
\begin{equation}
\alpha_{1}\leq\frac{\sum\limits_{i=1}^{n}\underline Z_{i}-n\underline{\mu}}{\sqrt
{n}\underline{\sigma}}, \ \frac{\sum\limits_{i=1}^{n}\bar{Z}_{i}-n\bar{\mu}%
}{\sqrt{n}\bar{\sigma}}\leq\alpha_{2}\Leftrightarrow T\left(
\begin{array}
[c]{c}%
\frac{\sum\limits_{i=1}^{n}\underline Z_{i}-n\underline{\mu}}{\sqrt{n}\underline{\sigma}%
}\\
\frac{\sum\limits_{i=1}^{n}\bar{Z}_{i}-n\bar{\mu}}{\sqrt{n}\bar{\sigma}}%
\end{array}
\right)  \in C, \label{appberry}%
\end{equation}
for some convex $C\subset\mathbb{R}^{2}$.

According to the multidimensional Berry-Esseen Theorem \cite{Da}, we get
\begin{equation}
\Bigg |P_{\nu}^{\infty}\Bigg (T\left(
\begin{array}
[c]{c}%
\frac{\sum\limits_{i=1}^{n}\underline Z_{i}-n\underline{\mu}}{\sqrt{n}\underline{\sigma}%
}\\
\frac{\sum\limits_{i=1}^{n}\bar{Z}_{i}-n\bar{\mu}}{\sqrt{n}\bar{\sigma}}%
\end{array}
\right)  \in C\Bigg )-\mathrm{Pr}\Bigg(\left(
\begin{array}
[c]{c}%
\hat B_1\\
\hat B_2%
\end{array}
\right)  \in C\Bigg)\Bigg |\leq\frac{K}{\sqrt{n}} \label{berry},%
\end{equation}

for some constant K, where $\left(
\begin{array}
[c]{c}%
\hat B_1\\
\hat B_2%
\end{array}
\right)  $ is the standard bivariate normal.

Define
\[
\left(
\begin{array}
[c]{c}%
B_1\\
B_2%
\end{array}
\right)  =T^{-1}\left(
\begin{array}
[c]{c}%
\hat B_1\\
\hat B_2%
\end{array}
\right)  .
\]

By (\ref{appberry}), we have
\begin{equation}
\mathrm{Pr}\Bigg(\left(
\begin{array}
[c]{c}%
\hat B_1\\
\hat B_2%
\end{array}
\right)  \in C\Bigg)=\mathrm{Pr}(\alpha_{1}\leq B_1,B_2%
\leq\alpha_{2})\label{normal}\\
=\mathrm{Pr}(-B_1\leq-\alpha_{1},B_2\leq\alpha_{2}%
)=\mathbf{N}_{2}(-\alpha_{1},\alpha_{2};-\rho).
\end{equation}

Then Theorem \ref{twosideclt} follows from (\ref{Z}), (\ref{appberry}),
(\ref{berry}), (\ref{normal}). This completes the proof.$\ \ \ \ \ \ \Box$

\begin{remark}
From equation (\ref{rho}), $\rho$ depends on  $M$. Actually,  $\underline\mu,\ \bar\mu,\ \underline\sigma,\ \bar\sigma$ depend on $M$ too. Recall that  $M$ is the least upper bound of $\{Y_{i}\}_{i=1}^\infty$, that is, $M=inf\{C \geq 0:\nu(|Y_i|\leq C)=1\}$. So $M$ is a parameter reflecting the distributional property of $\{Y_i\}_{i=1}^\infty$. But we claim that if $M$ is replaced by any number large than $M$ in equation (\ref{rho}), the value of  $\rho$ does not change. For simplicity, we suppose $Y_{i}\geq
0, \ i\geq1$. In this case,
\[
E^{P_{\nu}^{\infty}}(\underline Z_{i}\bar{Z}_{i}%
)=M\underline{\mu}-\int_{0}^{M}\int%
_{0}^{t_{2}}\nu(t_{1}\leq Y_{i}\leq t_{2})\mathrm{d}t_{1}\mathrm{d}t_{2}
\]

It's enough to prove
\begin{equation}
\label{rhoM}
(M+1)\underline{\mu}-\int_{0}^{M+1}\int_{0}^{t_{2}}\nu(t_{1}\leq Y_{i}\leq
t_{2})\mathrm{d}t_{1}\mathrm{d}t_{2}=M\underline{\mu}-\int_{0}^{M}\int_{0}^{t_{2}}\nu(t_{1}\leq Y_{i}\leq
t_{2})\mathrm{d}t_{1}\mathrm{d}t_{2}.
\end{equation}

And equation (\ref{rhoM}) follows from
\begin{align*}
& \ \ \ \ \ (M+1)\underline{\mu}-\int_{0}^{M+1}\int_{0}^{t_{2}}\nu(t_{1}\leq Y_{i}\leq
t_{2})\mathrm{d}t_{1}\mathrm{d}t_{2}\\
&  =(M+1)\underline{\mu}-\int_{0}^{M}\int_{0}^{t_{2}}\nu(t_{1}\leq Y_{i}\leq
t_{2})\mathrm{d}t_{1}\mathrm{d}t_{2}-\int_{M}^{M+1}\int_{0}^{t_{2}}\nu
(t_{1}\leq Y_{i}\leq t_{2})\mathrm{d}t_{1}\mathrm{d}t_{2}\\
&  =(M+1)\underline{\mu}-\int_{0}^{M}\int_{0}^{t_{2}}\nu(t_{1}\leq Y_{i}\leq
t_{2})\mathrm{d}t_{1}\mathrm{d}t_{2}-\int_{M}^{M+1}\int_{0}^{M}\nu(t_{1}\leq
Y_{i})\mathrm{d}t_{1}\mathrm{d}t_{2}\\
&  =(M+1)\underline{\mu}-\int_{0}^{M}\int_{0}^{t_{2}}\nu(t_{1}\leq Y_{i}\leq
t_{2})\mathrm{d}t_{1}\mathrm{d}t_{2}-\underline{\mu}\\
&  =M\underline{\mu}-\int_{0}^{M}\int_{0}^{t_{2}}\nu(t_{1}\leq Y_{i}\leq
t_{2})\mathrm{d}t_{1}\mathrm{d}t_{2}.
\end{align*}

\end{remark}

\begin{remark}
If $X_i(\omega_1, \omega_2,...)=I_{\{\omega_i\in A_1\}}, A_1\in \mathcal{B}(\Omega)$, then
\[
\underline\mu=\nu(A_1), \ \bar\mu=\nu^*(A_1), \ \underline\sigma=\nu(A_1)(1-\nu(A_1)), \  \bar\sigma=\nu^*(A_1)(1-\nu^*(A_1)),
\]
\begin{align*}
\rho=corr^{P_{\nu}^{\infty}}(\underline Z_{i},\bar{Z}_{i})&=\frac
{M^{2}-M\bar{\mu}+M\underline{\mu}-\rho^{\prime}-\underline{\mu}\bar{\mu}%
}{\underline{\sigma}\bar{\sigma}}\\
&=\frac{1-\nu^*(A_1)+\nu(A_1)-\int_{-1}^1\int_{-1}^{t_2}\nu(t_1\leq Y_i\leq t_2)\mathrm{d}t_1\mathrm{d}t_2-\nu(A_1)\nu^*(A_1)}{\sqrt{\nu(A_1)(1-\nu(A_1)\nu*(A_1)(1-\nu^*(A_1))}}\\
&=\frac{1-\nu^*(A_1)+\nu(A_1)-\int_{0}^1\int_{-1}^{0}\nu( Y_i\leq t_2)\mathrm{d}t_1\mathrm{d}t_2-\nu(A_1)\nu^*(A_1)}{\sqrt{\nu(A_1)(1-\nu(A_1)\nu^*(A_1)(1-\nu^*(A_1))}}\\
&=\frac{1-\nu^*(A_1)+\nu(A_1)-\nu(A_1^c)-\nu(A_1)\nu^*(A_1)}{\sqrt{\nu(A_1)(1-\nu(A_1)\nu^*(A_1)(1-\nu^*(A_1))}}\\
&=\frac{\nu(A_1)(1-\nu^*(A_1))}{\sqrt{\nu(A_1)(1-\nu(A_1)\nu^*(A_1)(1-\nu^*(A_1))}}\\
&=\frac{\nu(A_1)\nu(A_1^c)}{\sqrt{\nu(A_1)(1-\nu(A_1)\nu^*(A_1)(1-\nu^*(A_1))}}.
\end{align*}

By Theorem \ref{twosideclt},
\begin{align*}
&\ \ \ \ \ \lim_{n\rightarrow\infty}\nu^{\infty}(\alpha_{1}\sqrt{n}\underline{\sigma}+n\underline{\mu}%
\leq\sum\limits_{i=1}^{n}X_{i}\leq\alpha_{2}\sqrt{n}\bar{\sigma}+n\bar{\mu
})\\
&=\mathbf{N_{2}}(-\alpha_{1},\alpha_{2};-\rho)\\
&=P\Big(\left(
\begin{array}
[c]{c}%
B_1\\
B_2%
\end{array}
\right)\leq\left(
\begin{array}
[c]{c}%
-\alpha_1\\
\alpha_2%
\end{array}
\right)\Big)\\
&=P\Big(\left(
\begin{array}
[c]{c}%
\sqrt{\nu(A_1)(1-\nu(A_1))}B_1\\
\sqrt{\nu^*(A_1)(1-\nu^*(A_1))}B_2%
\end{array}
\right)\leq\left(
\begin{array}
[c]{c}%
-\sqrt{\nu(A_1)(1-\nu(A_1))}\alpha_1\\
\sqrt{\nu^*(A_1)(1-\nu^*(A_1))}\alpha_2%
\end{array}
\right)\Big),\\
\end{align*}
where $\left(
\begin{array}
[c]{c}%
B_1\\
B_2%
\end{array}
\right)$ is the  bivariate normal with zero mean and covariance matrix $\left(
\begin{array}
[c]{cc}%
1 & -\rho\\
-\rho & 1%
\end{array}
\right).  $ \\
Notice that $\left(
\begin{array}
[c]{c}%
\sqrt{\nu(A_1)(1-\nu(A_1))}B_1\\
\sqrt{\nu^*(A_1)(1-\nu^*(A_1))}B_2%
\end{array}
\right)$ is a bivariate normal with zero mean and covariance matrix  \\
$\Lambda=\left(
\begin{array}
[c]{cc}%
\nu(A_1)(1-\nu(A_1)) & -\nu(A_1)\nu(A_2)\\
-\nu(A_1)\nu(A_2) & \nu(A_2)(1-\nu(A_2))%
\end{array}
\right)  $
which coincident with (\ref{cltinter}).
\end{remark}

\begin{remark}
When $\nu$ is additive, Theorem \ref{twosideclt} degenerates into the
classical central limit theorem,
\[
\lim_{n\rightarrow\infty}\nu^{\infty}(\alpha_{1}<\frac{\sum\limits_{i=1}%
^{n}X_{i}-n\underline{\mu}}{\sqrt{n}\underline{\sigma}}\leq\alpha
_{2})=\mathrm{Pr}(\alpha_{1}<B\leq\alpha_{2}),
\]
where $B$ is the standard normal.
In fact, in this case
\[
E^{P_\nu^\infty}X_i=\underline\mu=\bar\mu, \ \ \ Var X_i=\underline\sigma=\bar\sigma,
\]
 and
\begin{align*}
E^{P_{\nu}^{\infty}}(\underline Z_{i}\bar{Z}_{i})&=M^{2}-M\bar{\mu}+M\underline{\mu}-\int_{-M}^{M}\int_{-M}^{t_{2}}\nu^\infty
(t_{1}\leq X_{i}\leq t_{2})\mathrm{d}t_{1}\mathrm{d}t_{2}\\
&=M^{2}-\int_{-M}^{M}\int_{-M}^{t_{2}}[\nu^\infty (X_{i}\leq t_{2})-\nu^\infty (X_{i}\leq t_{1})]\mathrm{d}t_{1}\mathrm{d}t_{2}\\
&=M^{2}-\int_{-M}^{M}\int_{-M}^{t_{2}}\nu^\infty (X_{i}\leq t_{2})\mathrm{d}t_{1}\mathrm{d}t_{2}+\int_{-M}^{M}\int_{-M}^{t_{2}}\nu^\infty (X_{i}\leq t_{1})\mathrm{d}t_{1}\mathrm{d}t_{2}\\
&=M^{2}-\int_{-M}^{M}\int_{-M}^{t_{2}}\nu^\infty (X_{i}\leq t_{2})\mathrm{d}t_{1}\mathrm{d}t_{2}+\int_{-M}^{M}\int^{M}_{t_{1}}\nu^\infty (X_{i}\leq t_{1})\mathrm{d}t_{2}\mathrm{d}t_{1}\\
&=M^{2}-\int_{-M}^{M}(t_2+M)\nu^\infty(X_i\leq t_2)\mathrm{d}t_2+\int_{-M}^{M}(M-t_1)\nu^\infty(X_i\leq t_1)\mathrm{d}t_1 \\
&=M^{2}-2\int_{-M}^{M}t\nu^\infty(X_i \leq t)\mathrm{d}t\\
&=M^{2}-\int_{-M}^{M}\nu^\infty(X_i \leq t)\mathrm{d}t^2\\
&=M^{2}-t^2\nu^\infty(X_i\leq t)\Big|_{-M}^{M}+\int_{-M}^{M}t^2\mathrm{d}\nu^\infty(X_i \leq t)\\
&=M^{2}-M^{2}+E^{P_\nu^\infty}(X_i^2)\\
&=E^{P_\nu^\infty}(X_i^2),\ \ i\geq 1.
\end{align*}
So $\rho=\frac{E^{P_{\nu}^{\infty}}(\underline Z_i\bar Z_i)-E^{P_{\nu}^{\infty}}(\underline Z_i)E^{P_{\nu}^{\infty}}(\bar Z_i)}{\sqrt{E^{P_{\nu}^{\infty}}(\underline Z_i^2)-(E^{P_{\nu}^{\infty}}(\underline Z_i))^2}\sqrt{E^{P_{\nu}^{\infty}}(\overline Z_i^2)-(E^{P_{\nu}^{\infty}}(\overline Z_i))^2}}=\frac{E^{P_\nu^\infty}(X_i^2)-\underline\mu\bar\mu}{\underline\sigma\bar\sigma}=\frac{Var X_i}{Var X_i}=1$ and our theorem comes back to the classical CLT.
\end{remark}

\renewcommand{\refname}


{\large References}

\bigskip

\begin{thebibliography}{99}                                                                                               %


\bibitem {Al}M. Allais, \emph{Le comportement de l'homme rationel devant le
risque: Critiqsuse des postulates et axiomes de l'ecole americaine},
Econometrica, 21(1953), pp.503-546.


\bibitem {CWL}Z. Chen, P. Wu and B. Li, \emph{A strong law of large numbers for
non-additive probabilities}, International Journal of Approximate Reasoning, 54(2013), pp.365-377.

\bibitem {Ch}G. Choquet, \emph{Theory of capacities},Ann.Inst.Fouier(Grenoble), 5(1954), pp.131-295.

\bibitem {Da}A. Dasgupta, \emph{Asymptotic Theory of Statistics and
Probability}, Springer(2008).

\bibitem {DW}J. Dow and S, Werlang, \emph{Laws of large numbers for non-additive
probabilities}, Economics Working Papers(Ensaios Economicos da EPGE), 226(1993).

\bibitem {El}D. Ellsberg, \emph{Risk, ambiguity, and the Savage axioms}, The quarterly journal of economics, 75(1961), pp. 643-669.

\bibitem {EJ1}L. Epstein and S.L.Ji, \emph{Ambiguous volatility and asset pricing in continuous time}, The Review of Financial Study, 26(2013), pp.1740-1786.

\bibitem {EJ2}L. Epstein and S.L.Ji, \emph{Ambiguous volatility, possibility and utility in continuous time}, Journal of Mathematical Economics, 50(2014), pp.269-282.

\bibitem {ES1}L. Epstein, H. Kaido and K. Seo, \emph{Robust confidence regions for incomplete models}, Econometrica, 84(2016), pp.1799-1838.

\bibitem {ESc}L. Epstein and D. Schneider, \emph{IID: independently and
indistingguishably distributed}, Journal of Economic Theory, 113(2003), pp.32-50.

\bibitem {ES}L. Epstein and K. Seo, \emph{Exchangeable Capacities, Parameters and Incomplete Theories}, Journal of Economic Theory, 157(2015), pp.879-917.

\bibitem {HC}F. Hu, Z. Chen, \emph{General laws of large numbers under sublinear expectations}. Communications in Statistics-Theory and Methods, 45(2016), pp.4215-4229.

\bibitem {HS}P. Huber, V. Strassen, \emph{Minimax tests and the Neyman-Pearson lemma for capacities}. The Annals of Statistics, 1(1973), pp.251-263.

\bibitem {HZ}Z. Hu, L. Zhou, \emph{Multi-dimensional central limit theorems and laws of large numbers under sublinear expectations}. Acta Mathematica Sinica. English Series, 31(2015), pp. 305-318.

\bibitem{LS}M. Li, Y. Shi, \emph{A general central limit theorem under sublinear expectations}. Science China Mathematics, 53(2010), pp.1989-1994.

\bibitem {MM}F. Macceroni and M. Marinacci, \emph{A strong law of large numbers for
capacities}. The Annals of Probability, 33(2005), pp.1171-1178.

\bibitem {Ma}M, Marnacci, \emph{Limit laws for non-additive probabilities and
their frequentist interpretation}, Journal of Economic Theory, 84(1999), pp.145-195.

\bibitem {Pe3}S. Peng,  \emph{G-expectation,G-Brownian motion and related
stochastic calculus of It$\hat o$ type}, Stochastic Analysis and
Applications, Abel Symposia, Springer,Berlin Heidelberg,(2007), pp.541-567.

\bibitem {Pe4}S. Peng, \emph{Muti-Dimensional G-Brownian motion and related
Stochastic Calculus under G-Expectation}, Stochastic Processes and their
Applications, 118(2008), pp.2223-2253.

\bibitem {Pe7}S. Peng, \emph{A new central limit theorem under sublinear expectations}. arxiv:0803.2656, (2008).

\bibitem {Pe8}S. Peng, \emph{Law of large numbers and central limit theorem under
nonlinear expectations}. arxiv:0702358, (2007).

\bibitem {Pe5}S. Peng, \emph{Survey on normal distributions, central limit
theorem, Brownian motion and the related stochastic calculus under sublinear
expectations}. Science in China Series A-Mathematics, 52(2009), pp.1391-1411.

\bibitem {Pe6}S. Peng, \emph{Nonlinear expectations and stochastic calculus
under uncertainty with robust central limit theorem and G-Brownian motion}.
arXiv:1002.4546, (2010).

\bibitem {PDJ}F. Philippe, G. Debs and J.Y. Jaffray, \emph{Decision making with
monotone lower probabilities of infinite order}. Mathematics of Operations Research, 24(1999), pp.767-784.

\bibitem {Sc}D. Schmeidler, \emph{Subjective probability and expected utility
without additivity}. Econometrica, 57(1989), pp.571-587.


\bibitem {WF}P. Walley and T. Fine, \emph{Towards a frequentist theory of upper and
lower probability}. The Annals of Statistics, 10(1982), pp.741-761.

\bibitem {WKl}Z. Wang and G. Klir, \emph{Generalized measure theory}, Springer Science and Business Media(2010).

\bibitem {ZC}D. Zhang, Z. Chen,\emph{A weighted central limit theorem under sublinear expectations}. Communications in Statistics-Theory and Methods, 43(2014), pp.566-577.



\end{thebibliography}
\end{document}